\documentclass[12pt, a4paper]{article}
\usepackage{amssymb}
\usepackage{graphicx}
\usepackage{psfrag}
\usepackage[spanish,activeacute]{babel}
\usepackage[latin1]{inputenc}

\newtheorem{defi}{\bf Definition}[subsection]
\newtheorem{lemma}{\bf Lemma}[subsection]
\newtheorem{prop}{\bf Proposition}[subsection]
\newtheorem{cor}{\bf Corollary}[subsection]
\newtheorem{teo}{\bf Theorem}[subsection]
\newtheorem{obs}{\bf Remark}[subsection]
\def\bp{\noindent{\it {\bf Proof:\ }}}
\newcommand{\ep}{\hfill$\Box$}

\def\R{I\kern -0.37 em R}
\def\N{I\kern -0.37 em N}

\newcommand{\Z}{{\bf Z}}

\begin{document}

\title{Expansive and fixed point free homeomorphisms of the plane.}
\author{J. Groisman\thanks{Research partially supported by ICTP, Trieste, Italy.}}
\maketitle
\begin{center}
{\small IMERL, Facultad de Ingenier\'{\i}a, Universidad de la
Rep\'ublica, Montevideo Uruguay\\ (e-mail: jorgeg@fing.edu.uy)}
\end{center}

\begin{abstract}
The aim of this work is to describe the set of  fixed point free homeomorphisms of the plane under certain expansive conditions.
\end{abstract}

\section{Introduction.}\label{I}
In \cite{G}, necessary and sufficient conditions for a
homeomorphism of the plane with one fixed point to be
topologically conjugate to a linear hyperbolic automorphism was
proved. The discovery of a hypothesis about the behavior of
Lyapunov functions at infinity was essential for this purpose. In
this work we will describe the set of fixed point free homeomorphisms of the plane which admit a
Lyapunov metric function with the same
condition imposed in the previously cited article. These homeomorphisms
are expansive respect to the Lyapunov metric $U$, meaning $U:\R^{2}
\times \R^{2} \rightarrow \R $ continuous and positive (i.e. it is
equal to zero only on the diagonal) and $W=\Delta (\Delta U) $
positive with $\Delta U(x,y)=U(f(x),f(y))-U(x,y)$.
As it is well known, expansive homeomorphisms on compact surfaces were
classified by Lewowicz in \cite{L1} and Hiraide in \cite{H}. As a matter of fact, we began
by studying whether some of the results obtained in the previously
cited article could be adapted to our new context (i.e. without
working in a compact environment but having the local compactness
of the plane). These arguments will allow us to construct singular
transverse foliations $\mathcal{F}^{s}$ and $\mathcal{F}^{u}$. The main result of this
article, theorem \ref{TP}, states: {\it Let $f$ be a homeomorphism of $\R^{2}$ which is fixed
point free and admits a Lyapunov function
$U:\R^{2}\times\R^{2}\rightarrow\R$ such that the following
properties hold:
\begin{itemize}
\item $U$ is a metric in $\R^{2} $ and induces the same topology
in the plane as the usual metric, \item For each point $x\in
\R^{2}$ and for each $k>0$ there exist points $y$ and $z$ on the
boundary of $B_{k}(x)$ such that $V(x,y)=U(f(x),f(y))-U(x,y)>0$ and
$V(x,z)=U(f(x),f(z))-U(x,z)<0$, respectively, \item $f$ does not
admit singularities,
\end{itemize}
$f$ is topologically conjugate to a translation of the plane if and only if $U$ admits condition {\bf HP} which states that given any compact
set $C\subset \R^{2} $ the following property holds: $$ \lim _{\|
x\| \rightarrow \infty } \frac{|V(x,y)-V(x,z)|}{W(x,y)} =0,$$
uniformly with $y,z$ in $C$.}
Let´s take as an example a translation of the plane and show that admits a Lyapunov metric function with the condition {\bf HP}. For this aim we will adapt some arguments used by White in \cite{WW}. Let $T$ be a translation of the plane defined by $T(x,y)=(x+1,y)$. Fix $\lambda >1$ and consider the Riemannian structure defined by $$(u,v)_{\sigma }=\lambda^{-2x} (u,e^{s})(v,e^{s}) +  \lambda^{2x} (u,e^{u})(v,e^{u}),$$ where $e^{s}=(1,0)$ and $e^{u}=(0,1)$. The metric induced on the plane by this Riemannian metric
is defined by $$d_{\sigma} (p,q)=\inf_{\gamma} \int_{[0,1]} \| \gamma'(t) \|_{\sigma }dt,$$ where $\gamma $ is an arc joining points $p$ and $q$. The topology induced by $d_{\sigma }$ on the plane is the same that is induced by the usual metric. Let $\gamma =(\gamma_{1} ,\gamma_{2} ):[0,1]\rightarrow \R^{2}$ be a curve joining points $p$ and $q$ of $\R^{2}$. We can consider two pseudo metrics $D^{s}$ and $D^{u}$ defined by $$D^{s}(p,q)=\inf_{\gamma } \int_{[0,1]} \lambda ^{-\gamma_{1}(t)}|(\gamma' (t),e^{s})|dt,$$ and $$D^{u}(p,q)=\inf_{\gamma } \int_{[0,1]} \lambda ^{\gamma_{1}(t)}|(\gamma' (t),e^{u})|dt.$$ By a simple computation we can prove that $$D^{s}(T(p),T(q))=\lambda^{-1} D^{s}(p,q),$$ and $$D^{u}(T(p),T(q))=\lambda D^{u}(p,q).$$ Consider the Lyapunov metric $D$ of $\R^{2}$ defined by $D=D^{s}+D^{u}$. Let´s test some conditions required for $D$:
\begin{itemize}
\item[(I)] {\bf Signs for $\Delta (D)$.}  $$ \Delta D(x,y)=D(f(x),f(y))-D(x,y) = $$
$$ D_{s}(f(x),f(y))-D_{s}(x,y)+ D_{u}(f(x),f(y))-D_{u}(x,y)=$$
$$ (\lambda -1)D_{u}(x,y) -(1-1/\lambda )D_{s}(x,y).$$
For every point $x\in \R^{2} $ and for every $k>0$, there are
points $y$ in the boundary of $B_{k}(x)$ such that $D_{u}(x,y)=0$
(this is true because the stable set separates the plane).
Therefore, $\Delta D(x,y)<0$ as we wanted. A similar argument lets
us find points $z\in \R^{2}$ such that $\Delta
D(x,z)>0$.

\item[(II)] {\bf Property HP.} Let $V=\Delta D$ and $W=\Delta^{2}
D$. We´ll prove that given any compact set $C\subset \R^{2} $, the following
property holds
$$ \lim _{\| x\| \rightarrow \infty } \frac{|V(x,y)-V(x,z)|}{W(x,y)}
=0,$$ uniformly with $y,z$ in $C$.\\ Since $$ \Delta^{2}
D(x,y)= \Delta D(f(x),f(y))-\Delta D(x,y)=$$ $$ (\lambda
-1)^{2}D_{u}(x,y) + (1-1/\lambda )^{2}D_{s}(x,y),$$ we can
conclude that $\Delta^{2} D(x,y)$ tends to infinity when $x$ tends
to infinity. Now,
$$|\Delta D(x,y)-\Delta D(x,z)|\leq $$ $$(\lambda
-1)|D_{u}(x,y)-D_{u}(x,z)|+(1-1/\lambda
)|D_{s}(x,y)-D_{s}(x,z)|\leq $$ $$(\lambda -1)D_{u}(z,y)
+(1-1/\lambda )D_{s}(z,y).$$ Then $|\Delta D(x,y)-\Delta D(x,z)|$
is uniformly bounded when points $y$ and $z$ lie on a compact set.
Then property {\bf HP} holds.
\end{itemize}

\section{Constructing foliations.}\label{CF}
Let $f:\R^{2} \rightarrow \R^{2}$ be a homeomorphism of the plane
that admits a Lyapunov metric function $U$, meaning $U:\R^{2}
\times \R^{2} \rightarrow \R $ continuous and positive (i.e. it is
equal to zero only on the diagonal) and $W=\Delta (\Delta U) $
positive with $\Delta U(x,y)=U(f(x),f(y))-U(x,y)$. In this section
we will resume some results of Lewowicz in \cite{L1}, and Groisman
in \cite{G}. We will work with the topology induced by a Lyapunov
function $U$ and define the $k$-stable set in the following way:
$$S_{k}(x)=\{y\in \R^{2}:\ U(f^{n}(x),f^{n}(y))\leq k,\ n\geq 0\}.$$
Similar definition for the $k$-unstable set. Let $f$ be a
homeomorphism of the plane that admits a Lyapunov function
$U:\R^{2}\times\R^{2}\rightarrow\R$ such that the following
properties hold:
\begin{itemize}

\item[(1)] {\bf $U$ is a metric in $\R^{2} $ and induces the same
topology in the plane as the usual metric.} Observe that given any
Lyapunov function it is possible to obtain another Lyapunov
function that verifies all the properties of a metric except,
perhaps, for the triangular property.

\item[(2)] {\bf Existence of both signs for the first difference
of $U$.} For each point $x\in \R^{2} $ and for each $k>0$ there
exist points $y$ and $z$ on the boundary of $B_{k}(x)$ such that
$V(x,y)=U(f(x),f(y))-U(x,y)>0$ and $V(x,z)=U(f(x),f(z))-U(x,z)<0$,
respectively.

\end{itemize}

\begin{obs}
A homeomorphism $f$ that admits a Lyapunov function $U$ defined at
$\R^{2}\times\R^{2}$ is U-expansive. This means that given two
different points of the plane $x,y$ and given any $k>0$, there
exists $n\in \Z $ such that $$U(f^{n}(x),f^{n}(y))>k.$$
\end{obs}

\bp Let $x$ and $y$ be two different points of the plane such that
$V(x,y)=U(f(x),f(y))-U(x,y)>0$. Since $\Delta V>0$, then
$V(f^{n}(x),f^{n}(y))>V(x,y)$ holds for $n>0$. This means that
$U(f^{n}(x),f^{n}(y))$ grows to infinity, since
$$U(f^{n}(x),f^{n}(y))=U(x,y)+\sum_{j=0}^{n-1} V(f^{j}(x),f^{j}(y))>$$
$$U(x,y)+nV(x,y).$$ Thus, given $k>0$ there exists $n\in \N $ such that
$$U(f^{n}(x),f^{n}(y))>k.$$ By using similar arguments we can prove the case
when $V(x,y)=U(f(x),f(y))-U(x,y)<0$. If $V(x,y)=0$, then
$V(f(x),f(y))>0$ and this is precisely our first case.\ep

\begin{defi}
Let $f:\R^{2} \rightarrow \R^{2}$ be a homeomorphism of the plane
that admits a Lyapunov metric function $U$. A point $x\in \R^{2} $
is a stable (unstable) point if given any $k'>0$ there exists
$k>0$ such that for every $y\in B_{k}(x)$, it follows that
$U(f^{n}(x),f^{n}(y))<k'$ for each $n\geq 0$ ($n\leq 0$).
\end{defi}

\begin{obs}
Property $(2)$ for $U$ implies the non-existence of stable
(unstable) points.
\end{obs}

\bp Given the existence of both signs for
$V(x,y)=U(f(x),f(y))-U(x,y)$ in any neighborhood of $x$, we can
state that for each $k>0$, there exists a point $y$ in $B_{k}(x)$
such that $V(x,y)>0$. Since $\Delta V>0$, we can state that
$V(f^{n}(x),f^{n}(y))>V(x,y)$ for $n>0$, so $U(f^{n}(x),f^{n}(y))$
grows to infinity. Thus, there are no stable points. We can use
similar arguments for the unstable case.\ep

\begin{lemma}
Let $A$ be an open set of $\R^{2} $ with $x\in A\subset B_{k}(x)$.
There exists a compact connected set $C$ with $x\in C\subset
\overline{A} $, $C\cap \partial (A)\neq \emptyset $ such that, for
all $y\in C$ and $n\geq 0$, $U(f^{n}(x),f^{n}(y))\leq k $ holds.
\end{lemma}

\bp See section $3$ \cite{G}

For $x\in \R^{2}$ and $k>0$, let $S_{k}(x)$ be the $k -$stable set
for $x$, defined by $$S_{k}(x)=\{ y\in \R^{2} :\
U(f^{n}(x),f^{n}(y))\leq k ,\ n\geq 0\} .$$

\begin{lemma}
Let us consider $0<k'<k $. There exists $\sigma >0$ such that if
$y\in S_{k}(x)$ and $U(x,y)<\sigma $, then $y\in S_{k'}(x)$.
\end{lemma}

\bp See section $3$ \cite{G} \\ Let $C_{k}(x)$ be the connected
component of $S_{k}(x)$ that contains $x$.
\begin{lemma}
$C_{k}(x)$ is locally connected at $x$.
\end{lemma}
\bp (See Lemma $2.3$, \cite{L1})
\begin{cor}
For each $x$ in $\R^{2} $, $C_{k }(x)$ is connected and locally
connected.
\end{cor}
\bp (See \cite{L1})
\begin{cor}
For any $x$ in $\R^{2} $ and any pair of points $p$ and $q$ in
$C_{k}(x)$, there exists an arc included in $C_{k}(x)$ that joins
$p$ and $q$.
\end{cor}
\bp (See Topology, Kuratowski, \cite{K}, section $50$)
\begin{prop}
Except for a discrete set of points, that we shall call singular,
every $x$ in $\R^{2} $ has local product structure. The stable
(unstable) sets of a singular point $y$ consists of the union of
$r$ arcs, with $r\geq 3$ that meet only at $y$. The stable
(unstable) arcs separate unstable (stable) sectors.
\end{prop}
\bp (See Section $3$, \cite{L1})
\begin{obs}
The neighborhood's size where there exists a local product
structure may become arbitrarily small. However we are able to
extend these stable and unstable arcs getting curves that we will
denote as $W^{s}(x)$ and $W^{u}(x)$, respectively. If two points
$y$ and $z$ belong to $W^{s}(x)$ ($W^{u}(x)$), then
$U(f^{n}(y),f^{n}(z)<k$ for some $k>0$ and for all $n\geq 0$
($n\leq 0$).
\end{obs}
The following lemmas refer to these stable and unstable curves.
\begin{lemma} \label{TT}
Let $f$ be a homeomorphism of the plane which verifies the
conditions of this section. Then stable and unstable curves
intersect each other at most once.
\end{lemma}
\bp If they intersect each other more than once, we would
contradict $U$-expansiveness: if two different points $x$ and $y$
belong to the intersection of a stable and an unstable curve, then
there exists $k_{0}>0$ such that $U(f^{n}(x),f^{n}(y)<k_{0}$ for
all $n\in \Z $.\ep
\begin{lemma}
Every stable (unstable) curve separates the plane.
\end{lemma}
\bp See section $3$ \cite{G}

Finally, we have conditions to state the following
theorem:
\begin{teo}
Let $f$ be a homeomorphism of the plane that admits a Lyapunov
function $U:\R^{2}\times\R^{2}\rightarrow\R$ such that the
following properties hold:
\begin{itemize}

\item $U$ is a metric in $\R^{2} $ and induces the same topology
in the plane as the usual metric.

\item For each point $x\in \R^{2} $ and for each $k>0$ there exist
points $y$ and $z$ on the boundary of $B_{k}(x)$ such that
$V(x,y)=U(f(x),f(y))-U(x,y)>0$ and $V(x,z)=U(f(x),f(z))-U(x,z)<0$,
respectively.

\end{itemize}
Then, $f$ admits transverse singular foliations $\mathcal{F}^{s}$
and $\mathcal{F}^{u}$. Leaves of $\mathcal{F}^{s}$
($\mathcal{F}^{u}$) are the stable (unstable) curves constructed in
this section.
\end{teo}

\section{Foliation description.}\label{DF}
In this section we will describe those foliations
introduced in the previous section with one additional condition
for the Lyapunov function $U$:\\ {\bf Property HP.} Let $V=\Delta
U$ and $W=\Delta^{2} U$. Given any compact set $C\subset \R^{2} $,
the following property holds
$$ \lim _{\| x\| \rightarrow \infty } \frac{|V(x,y)-V(x,z)|}{W(x,y)}
=0,$$ uniformly with $y,z$ in $C$.\\ In \cite{G} these foliations
were characterized in the case when the homeomorphism $f$ does not
admit singularities and it has a fixed point. Now we will
generalize these results for the case when $f$ is fixed point free.

\begin{teo} \label{TA}
Let $f$ be a homeomorphism of the plane such that the following
conditions hold:
\begin{itemize}
\item $f$ admits a Lyapunov metric function $U:\R^{2} \times
\R^{2} \rightarrow \R $. This metric induces in the plane the same
topology as the usual distance; \item $f$ has no singularities;
\item for each point $x\in \R^{2} $ and any $k>0$, there exist
points $y$ and $z$ in the boundary of $B_{k}(x)$ such that
$V(x,y)=U(f(x),f(y))-U(x,y)>0$ and $V(x,z)=U(f(x),f(z))-U(x,z)<0$,
where $B_{k}(x)=\{ y\in \R^{2} \ / U(x,y)\leq k\} $. \item {\bf
Property HP.} Given any compact set $C\subset \R^{2} $ the
following property holds: $$ \lim _{\| x\| \rightarrow \infty }
\frac{|V(x,y)-V(x,z)|}{W(x,y)} =0,$$ uniformly with $y,z$ in $C$.
\end{itemize}
Then, there exist transverse foliations $\mathcal{F}^{s}$ and
$\mathcal{F}^{u}$ such that every leaf of $\mathcal{F}^{s}$
intersects every leaf of $\mathcal{F}^{u}$.
\end{teo}
\bp Let $\mathcal{F}^{s}$ and $\mathcal{F}^{u}$ be the stable and
unstable foliations constructed in the last section. Let $x$ an
arbitrarily point of $\R^{2}$ and denote by $W^{s}(x)$ and
$W^{u}(x)$ the leaves of $\mathcal{F}^{s}$ and $\mathcal{F}^{u}$
through $x$. We will divide this proof in two steps:
\begin{itemize}
\item {\it If $W^{u}(y)\cap W^{s}(x)\neq \emptyset $ then
$W^{s}(y)\cap W^{u}(x)\neq \emptyset $, for every point $y$ in
$\R^{2}$.} If not let $z$ be the first point of $W^{u}(y)$ such
that $W^{s}(z)\cap W^{u}(x)=\emptyset $ (see fig. \ref{A}).

\begin{figure}[htb]\label{A}
\psfrag{x}{$x$} \psfrag{Wux}{$W^{u}(x)$} \psfrag{Wsx}{$W^{s}(x)$}
\psfrag{z}{$z$} \psfrag{Wuz}{$W^{u}(z)$} \psfrag{Wsz}{$W^{s}(z)$}
\begin{center}
\includegraphics[scale=0.2]{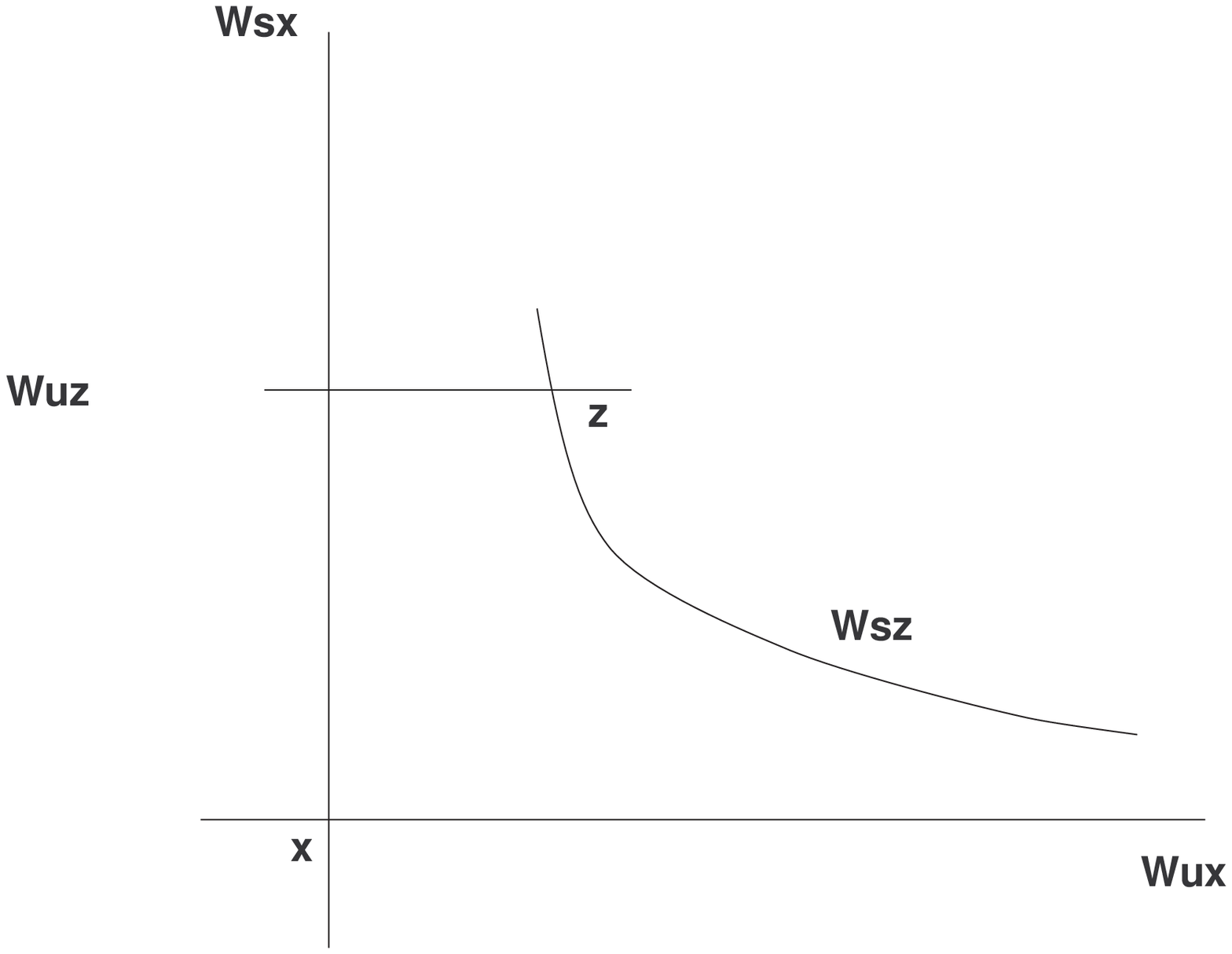}
\caption{\label{A}}
\end{center}
\end{figure}

Let $(z_{n})$ be a sequence of $W^{u}(y)$ such that $z_{n}$
converges to $z$ and $W^{s}(z_{n})\cap W^{u}(x)\neq \emptyset $,
for all $n$. As $n$ grows, the behavior of $W^{s}(z_{n})$ has two
possibilities as shown in fig \ref{B}.

\begin{figure}[htb]
\psfrag{x}{$x$} \psfrag{z}{$z$} \psfrag{zn}{$z_{n}$}
\psfrag{qn}{$q_{n}$} \psfrag{wn}{$w_{n}$}
\begin{center}
\includegraphics[scale=0.3]{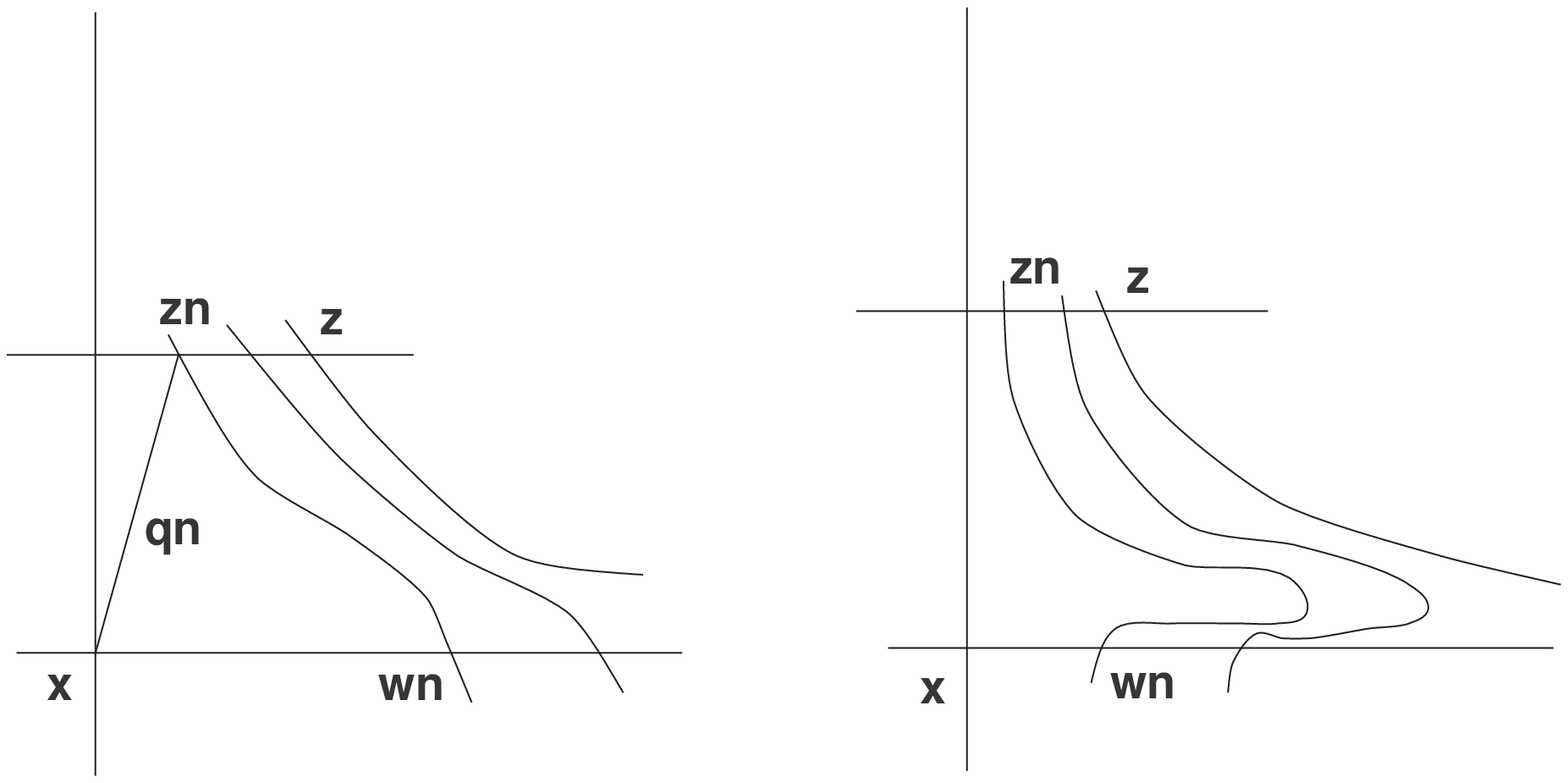}
\caption{\label{B}}
\end{center}
\end{figure}

In the first situation we could find points
$w_{n}=W^{s}(z_{n})\cap W^{u}(x)$ arbitrarily close to infinity.
Since for each $n$ we have that $V(w_{n},z_{n})<0$ (because they
belong to the same stable leaf) and $V(w_{n},x)>0$ (because they
belong to the same unstable leaf), there exists a point $q_{n}$,
belong to the segment line $z_{n}x$, such that $V(w_{n},q_{n})=0$.
Since $z_{n}$ and $q_{n}$ belong to a compact set for all $n$, we
can apply condition HP and find $n\in \N$ such that
$$\frac{|V(w_{n},z_{n})-V(w_{n},q_{n})|}{W(w_{n},z_{n})} <1,$$ which
implies that $$W(w_{n},z_{n})+ V(w_{n},z_{n})>0, $$ and then
$$V(f(w_{n}),f(z_{n}))>0.$$ This yields a contradiction since points
$f(w_{n}),f(z_{n})$ are in the same stable leaf.\\ In the second
situation (see fig \ref{C}), the set $\{w_{n}: n\in \N \}$ is
bounded. Let us consider, as figure \ref{C} shows, the segment
line $w_{n}z_{n}$ and let $h_{n}$ be the point of the arc belong to $W^{s}(z_{n})$ determined by the points $w_{n}$ and $z_{n}$, which is farthest from the point $x$.

\begin{figure}[htb]
\psfrag{x}{$x$} \psfrag{z}{$z$} \psfrag{zn}{$z_{n}$}
\psfrag{qn}{$q_{n}$} \psfrag{wn}{$w_{n}$} \psfrag{hn}{$h_{n}$}
\psfrag{rn}{$r_{n}$} \psfrag{Wuhn}{$W^{u}(h_{n})$}
\begin{center}
\includegraphics[scale=0.3]{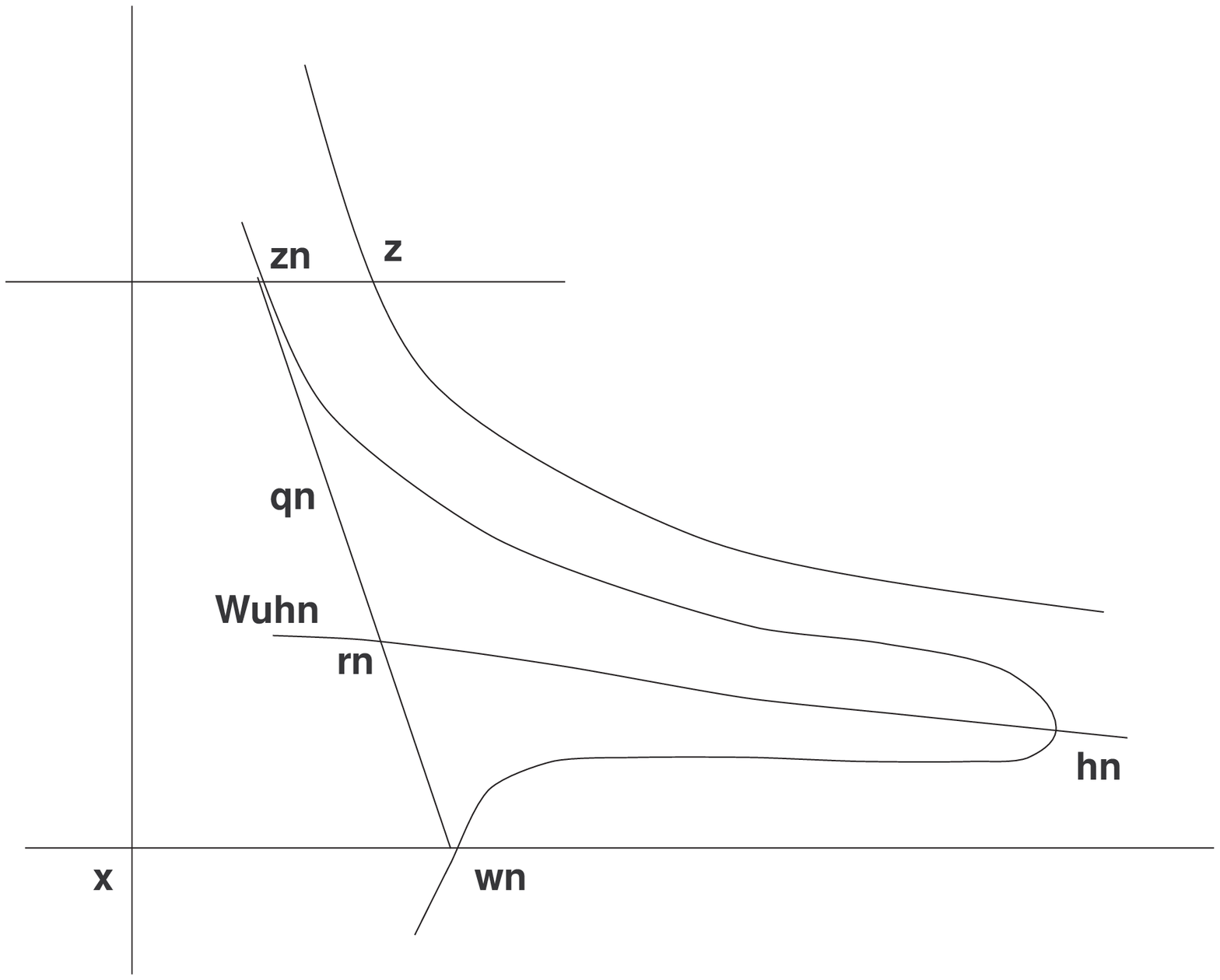}
\caption{\label{C}}
\end{center}
\end{figure}

$W^{u}(h_{n})$ must intersect segment $w_{n}z_{n}$. Otherwise, it
would cut $W^{s}(z_{n})$ more than once. Let $r_{n}$ be that
intersection point. We want to apply our condition {\bf HP}.
Observe that points $r_{n},z_{n}$ would be in a compact set for
all $n\in \N$ and $h_{n}$ tends to infinity when $z_{n}$ tends to
$z$. $V(h_{n},r_{n})>0$ because they belong to the same unstable
leaf, and $V(h_{n},z_{n})<0$ because they belong to the same
stable leaf. So, there exists a point $q_{n}$ that belongs to
segment line $z_{n}r_{n}$ such that $V(h_{n},q_{n})=0$. Then $$
\lim _{ n \rightarrow \infty }
\frac{|V(h_{n},z_{n})-V(h_{n},q_{n})|}{W(h_{n},z_{n})} =0.$$
Therefore, we can choose $h_{n}$ such that
$$W(h_{n},z_{n})+V(h_{n},z_{n})>0,$$ which implies that
$$V(f(h_{n}),f(z_{n}))>0.$$ This contradicts the fact that points
$f(h_{n}),f(z_{n})$ are in the same stable leaf.

\item {\it $W^{u}(y)\cap W^{s}(x)\neq \emptyset $ and
$W^{s}(y)\cap W^{u}(x)\neq \emptyset $, for every point $y$ in
$\R^{2}$.} Let us consider the set $A$ consisting of the points whose
stable (unstable) leaf intersects the unstable (stable) leaf of
point $x$. It is clear that $A$ is open. Let us prove that it is
also closed. Let $(q_{n})$ be a sequence of $A$, convergent to
some point $q$ (see figure \ref{fig11}).

\begin{figure}[htb]
\psfrag{p}{$p$} \psfrag{Wuq}{$W^{u}(q)$} \psfrag{Wsq}{$W^{s}(q)$}
\psfrag{Wuqn0}{$W^{u}(q_{n_{0}})$}
\psfrag{Wsqn0}{$W^{s}(q_{n_{0}})$} \psfrag{V(q)}{$V(q)$}
\psfrag{q}{$q$} \psfrag{qn0}{$q_{n_{0}}$}
\psfrag{alphan0}{$\alpha_{n_{0}}$}
\begin{center}
\includegraphics[scale=0.25]{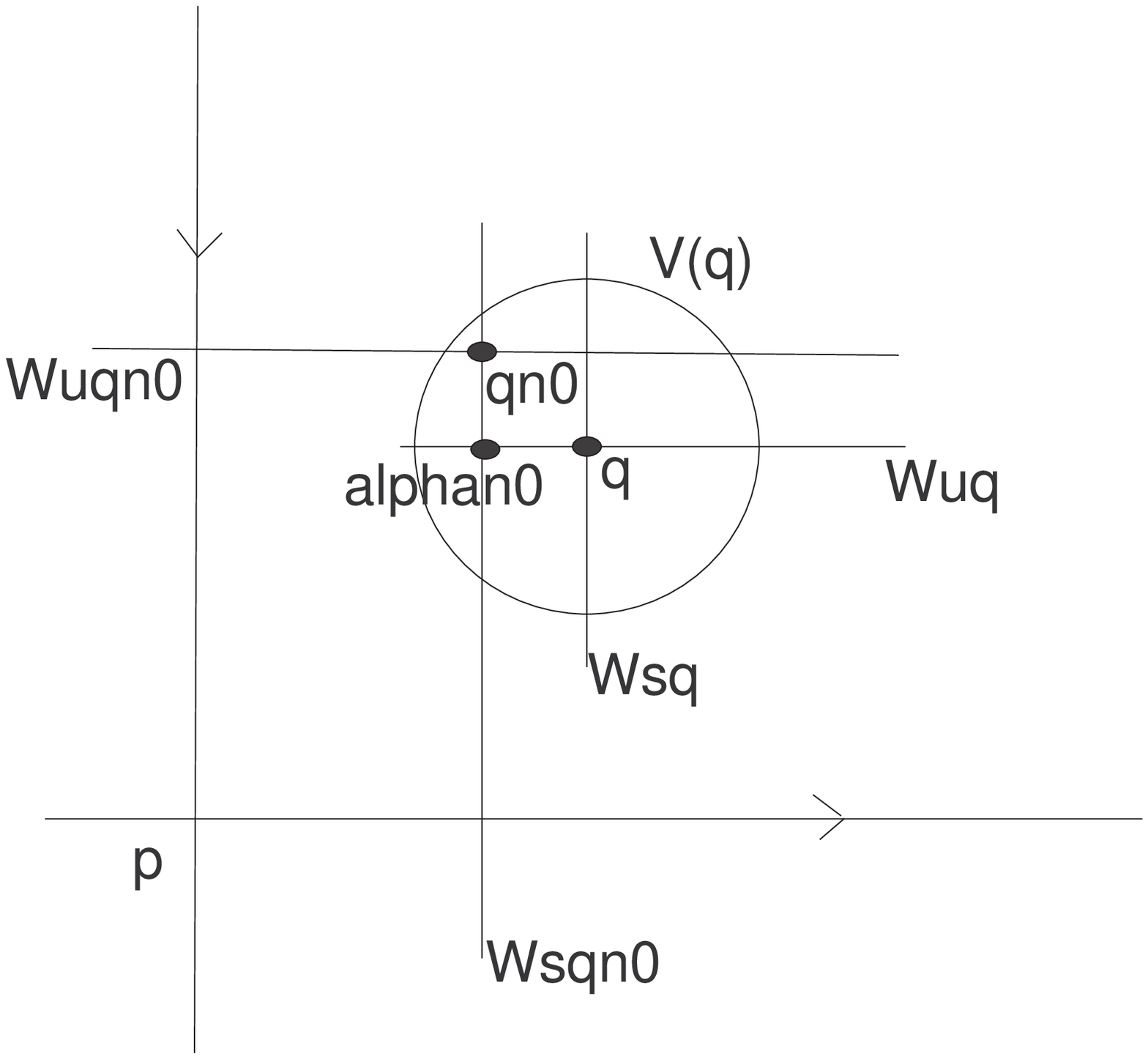}
\caption{\label{fig11}Coordinates}
\end{center}
\end{figure}

Let $V(q)$ be a neighborhood of $q$ with local product structure.
Let us consider $q_{n_{0}}\in V(q)$. So, we have that
$W^{s}(q_{n_{0}})\cap W^{u}(q)=\alpha_{n_{0}} $ as a consequence
of the local product structure and $W^{s}(q_{n_{0}})\cap
W^{u}(p)\neq \emptyset $ since $q_{n_{0}}\in A$. But then
$\alpha_{n_{0}} $ is a point in $W^{s}(q_{n_{0}})$ that cuts the
unstable leaf of point $x$, and then, applying the previous step
we have that $W^{u}(q)=W^{u}(\alpha_{n_{0}} )$ must cut the stable
leaf of point $x$. A similar argument lets us prove that the
stable leaf of $q$ must cut the unstable leaf of point $x$.
Therefore $q$ belongs to the set $A$ and consequently $A$ is closed.
Then $A$ is the whole plane.
\end{itemize}
Since $x$ is an arbitrary point, this proof is finished.\ep

\section{Main section.}\label{EW}
Let $f$ be a fixed point free
homeomorphism of $\R^{2}$. Brouwer's translation theorem (see \cite{B}, \cite{F}) asserts that if $f$ preserves orientation, then every $x_{0}\in \R^{2}$ is contained in a domain of
translation for $f$, i.e. an open connected subset of $\R^{2}$
whose boundary is $L\cup f(L)$ where $L$ is the image of a proper
embedding of $\R$ in $\R^{2}$, such that $L$ separates $f(L)$ and
$f^{-1}(L)$. The purpose of this section is to find some domain of translation of a fixed point free homeomorphism of $\R^{2}$ which admits a Lyapunov metric function $U$ such that properties presented in the last section hold. This situation will allow us to prove that the homeomorphisms discussed in theorem \ref{TA} are topologically conjugate to a translation of the plane.

\begin{lemma} \label{LP}
Let $f$ be a homeomorphism of $\R^{2}$ which admits a Lyapunov function $U:\R^{2}\times\R^{2}\rightarrow\R$ such that the following
properties hold:
\begin{itemize}
\item $U$ is a metric in $\R^{2} $ and induces the same topology
in the plane as the usual metric; \item for each point $x\in
\R^{2}$ and for each $k>0$ there exist points $y$ and $z$ on the
boundary of $B_{k}(x)$ such that $V(x,y)=U(f(x),f(y))-U(x,y)>0$ and
$V(x,z)=U(f(x),f(z))-U(x,z)<0$, respectively.
\end{itemize}
Let $W^{s}$ be an arbitrary non-invariant leaf of
$\mathcal{F}^{s}$ such that $W^{s}$ separates $f(W^{s})$ and
$f^{-1}(W^{s})$. Define $D$ as the open connected subset of
$\R^{2}$ whose boundary is $W^{s}\cup f(W^{s})$ and the open set
$U=\bigcup f^{k}(\overline{D})$ with $k\in \Z$. If $U\neq \R^{2}$
then the boundary of $U$ consists of the union of leaves
of $\mathcal{F}^{s}$. Similarly for the unstable case.
\end{lemma}
\bp If $U\neq \R^{2}$ then there exists a point $p$ and a sequence
$(x_{n})$ such that $x_{n}\in f^{n}(W^{s})\  n>0(n<0) $ and $\lim
x_{n}=p$. Let $W^{s}(p)$ be the leaf of $\mathcal{F}^{s}$
through $p$. Using continuity of $\mathcal{F}^{s}$ respect to the
initial point we could state that every point of $W^{s}(p)$
is the limit of a sequence $(y_{n})$ such that $y_{n}\in
f^{n}(W^{s})\ n>0 (n<0)$. Thus, $W^{s}(p)$ is included in the
boundary of $U$. If $W^{s}(p)$ is not invariant then all the iterates of $W^{s}(p)$ are also included in the boundary of $U$.\ep

\begin{figure}[htb]
\psfrag{x}{$x$} \psfrag{Ws}{$W^{s}$} \psfrag{fWs}{$f(W^{s})$}
\psfrag{f-1Ws}{$f^{-1}(W^{s})$} \psfrag{Wux}{$W^{u}(x)$}
\begin{center}
\includegraphics[scale=0.3]{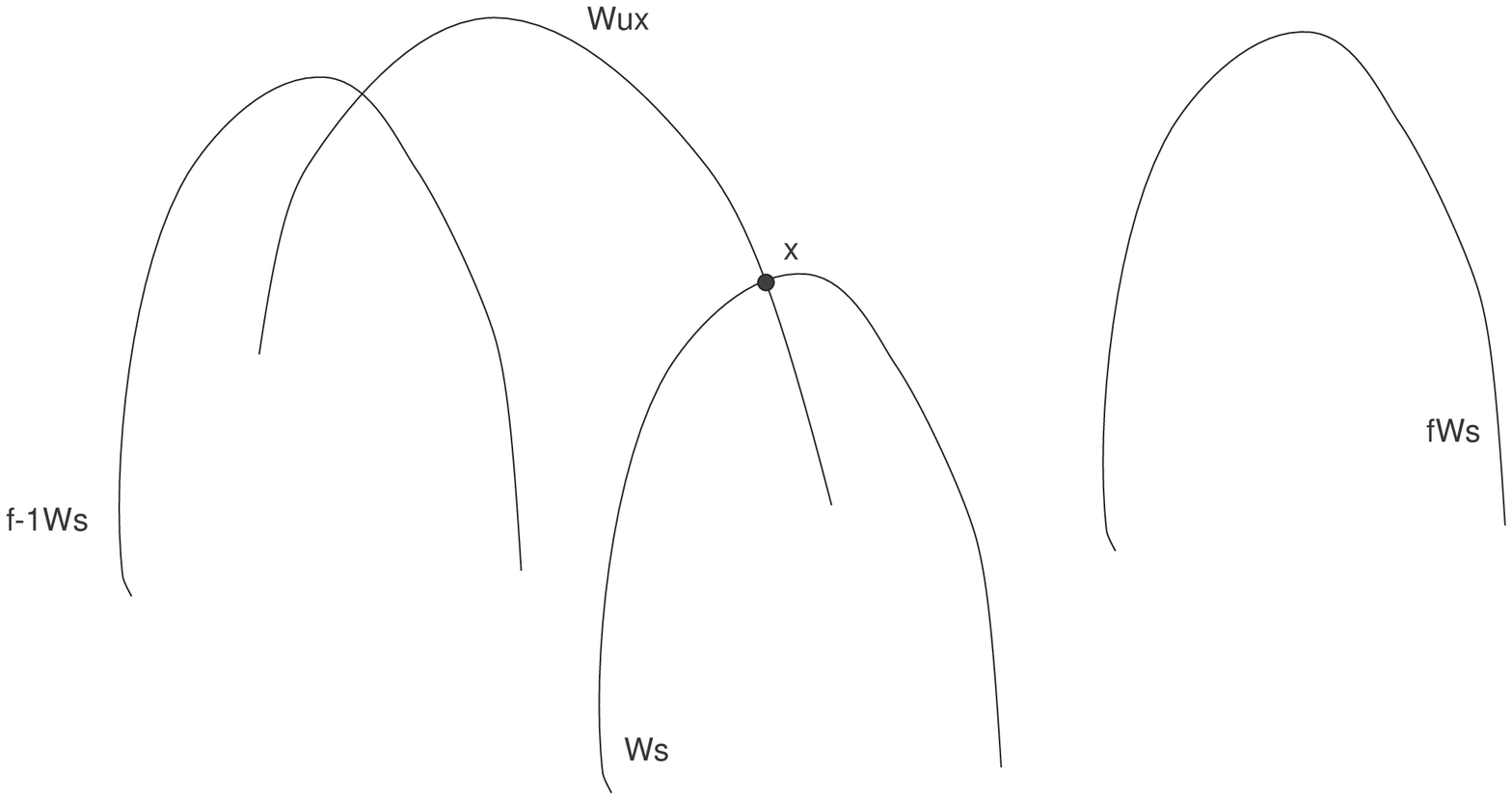}
\caption{\label{SEP}}
\end{center}
\end{figure}

\begin{lemma} \label{LP1}
Let $f$ be a homeomorphism of $\R^{2}$
which is fixed point free and admits a Lyapunov function
$U:\R^{2}\times\R^{2}\rightarrow\R$ such that the following
properties hold:
\begin{itemize}
\item $U$ is a metric in $\R^{2} $ and induces the same topology
in the plane as the usual metric, \item for each point $x\in
\R^{2}$ and for each $k>0$ there exist points $y$ and $z$ on the
boundary of $B_{k}(x)$ such that $V(x,y)=U(f(x),f(y))-U(x,y)>0$ and
$V(x,z)=U(f(x),f(z))-U(x,z)<0$, respectively,
\item $U$ admits condition (HP), \item $f$ has no
singularities.
\end{itemize}
If $W^{s}$ ($W^{u}$) is an arbitrary non invariant leaf of
$\mathcal{F}^{s}$ ($\mathcal{F}^{u}$), then $W^{s}$ ($W^{u}$)
separates $f^{2}(W^{s})$ ($f^{2}(W^{u})$) and $f^{-2}(W^{s})$
($f^{-2}(W^{u})$).
\end{lemma}
\bp Let $W^{s}$ be an arbitrary non invariant leaf of
$\mathcal{F}^{s}$ and consider the leaves $f(W^{s})$ and
$f^{-1}(W^{s})$. Let's suppose that none of these three leaves
separates the others (see fig. \ref{SEP}). 
Take any point $x$ in $W^{s}$. Applying theorem \ref{TA} we
have that $W^{u}(x)$ must intersect transversally $f(W^{s})$ and
$f^{-1}(W^{s})$. Also, recall that in lemma \ref{TT} we proved
that stable and unstable leaves intersect each other at most once.
Then $W^{u}(x)$ goes from one component to the other determined by
$W^{s}$ only once (remember that $W^{s}$ separates the plane).
Since we are assuming that $W^{s}$ does not separate $f(W^{s})$
and $f^{-1}(W^{s})$, then these two leaves are in the same
component determined by $W^{s}$. So, if $W^{u}(x)$ intersects
$f(W^{s})$ then it can not intersect $f^{-1}(W^{s})$
(because if it does, we would have either more than one intersection
between $W^{u}(x)$ and $f(W^{s})$ or an auto intersection
of $W^{u}(x)$). That yields a contradiction. Then, one of the three stable leaves separates the other two. Then $W^{s}$ must separate $f^{2}(W^{s})$ and $f^{-2}(W^{s})$.\ep

\begin{figure}[htb] \psfrag{a}{$a$} \psfrag{x}{$x$} \psfrag{y}{$y$} \psfrag{Ws}{$W^{s}$} \psfrag{fWs}{$f(W^{s})$}
    \psfrag{f2Ws}{$f^{2}(W^{s})$} \psfrag{Wsz}{$W^{s}(z)$} \psfrag{x4}{$x_{4}$} \psfrag{x6}{$x_{6}$} \psfrag{f3}{$f^{3}$} \psfrag{f4}{$f^{4}$} \psfrag{f5}{$f^{5}$} \psfrag{f6}{$f^{6}$} \begin{center} \includegraphics[scale=0.25]{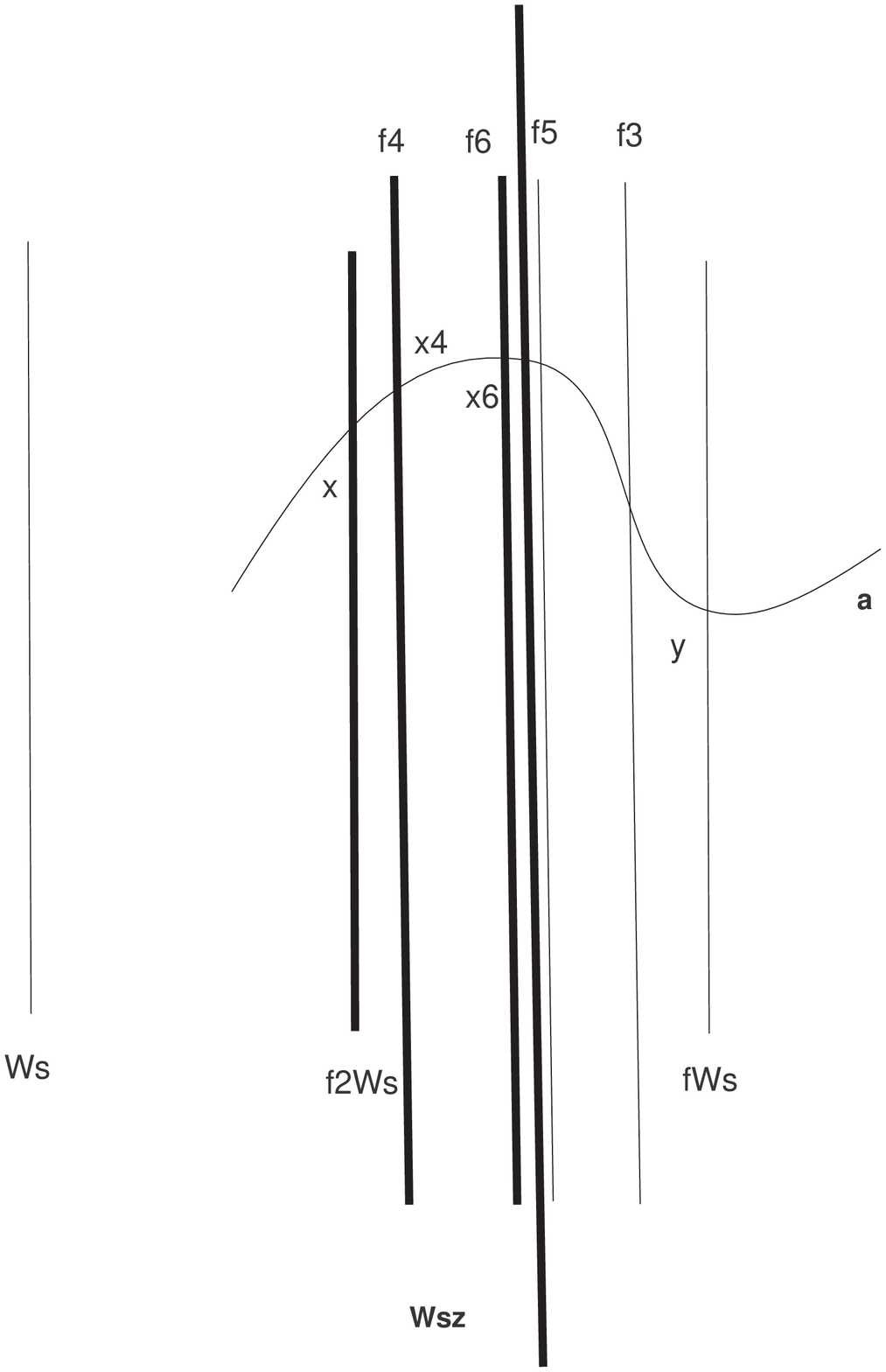} \caption{\label{H}} \end{center} \end{figure}

\begin{lemma} \label{LP2}
Let $f$ be a homeomorphism of $\R^{2}$
which is fixed point free and admits a Lyapunov function
$U:\R^{2}\times\R^{2}\rightarrow\R$ such that the following
properties hold:
\begin{itemize}
\item $U$ is a metric in $\R^{2} $ and induces the same topology
in the plane as the usual metric, \item for each point $x\in
\R^{2}$ and for each $k>0$ there exist points $y$ and $z$ on the
boundary of $B_{k}(x)$ such that $V(x,y)=U(f(x),f(y))-U(x,y)>0$ and
$V(x,z)=U(f(x),f(z))-U(x,z)<0$, respectively,
\item $U$ admits condition (HP), \item $f$ has no
singularities.
\end{itemize}
If $W^{s}$ ($W^{u}$) is an arbitrary non invariant leaf of
$\mathcal{F}^{s}$ ($\mathcal{F}^{u}$) such that $W^{s}$ ($W^{u}$)
does not separate $f(W^{s})$ ($f(W^{u})$) and $f^{-1}(W^{s})$
($f^{-1}(W^{u})$), then $f^{2n}(W^{s})$ or $f^{-2n}(W^{s})$, $n>0$, is a sequence of stable leaves converging to a $f^{2}$-invariant stable leaf.
\end{lemma}
\bp Let's suppose that $f(W^{s})$ separates $W^{s}$ and $f^{-1}(W^{s})$. As fig. \ref{H} shows let us consider an unstable arc $a$ joining a point $x$ of $f(W^{s})$ with a point $y$ of $W^{s}$.
    Since $f^{n}(W^{s})$ must separate $f^{n-1}(W^{s})$ and $f^{n-2}(W^{s})$ we conclude that $f^{n}(W^{s})$ must intersect the compact arc $a$ in a point $x_{n}$. Using lemma \ref{LP1}, we can state that sequence $x_{2n}$ is monotone and bounded. Let's denote by $z$ its limit. Since $W^{s}(z)$ separates $f^{-2}(W^{s})$ and $f^{2}(W^{s})$ we conclude that $W^{s}(z)$ is $f^{2}$-invariant. The other case is analogous.\ep
\begin{teo} \label{TP}
Let $f$ be a homeomorphism of $\R^{2}$
which is fixed point free and admits a Lyapunov function
$U:\R^{2}\times\R^{2}\rightarrow\R$ such that the following
properties hold:
\begin{itemize}
\item $U$ is a metric in $\R^{2} $ and induces the same topology
in the plane as the usual metric, \item for each point $x\in
\R^{2}$ and for each $k>0$ there exist points $y$ and $z$ on the
boundary of $B_{k}(x)$ such that $V(x,y)=U(f(x),f(y))-U(x,y)>0$ and
$V(x,z)=U(f(x),f(z))-U(x,z)<0$, respectively, \item $f$ has no
singularities.
\end{itemize}
Then, $f$ is topologically conjugate to a translation of the plane
if and only if $U$ admits condition {\bf HP}.
\end{teo}

\bp Let $W^{s}$ be an arbitrary non invariant leaf of
$\mathcal{F}^{s}$. Let's divide the proof in two cases:
\begin{itemize}
\item $W^{s}$ separates $f^{-1}(W^{s})$ and $f(W^{s})$. Define $D$ as the open connected subset of
$\R^{2}$ whose boundary is $W^{s}\cup f(W^{s})$ and the open set
$U=\bigcup f^{k}(\overline{D})$ with $k\in \Z$. $U$ is an open set invariant under $f$ such that the restriction of $f$ to $U$ is topologically conjugate to a translation of $\R^{2}$. If $U=\R^{2}$, then the theorem is proved. If not, applying lemmas \ref{LP} and \ref{LP1} there exists a $f^{2}$-invariant stable leaf $S$ in the boundary of $U$ (since $S$ separates $f^{2}(S)$ and $f^{-2}(S)$,
it is easy to prove that $S$ is $f^{2}$-invariant). Let's consider an arbitrary point $x$ of $S$ and its unstable leaf $W^{u}$. If $W^{u}$ separates $f^{-1}(W^{u})$ and $f(W^{u})$, consider the set $D'$ with boundary
$W^{u}(x)\cup f(W^{u}(x)$ and $U'=\bigcup f^{k}(\overline{D'})$ with $k\in \Z$. If $U'=\R^{2}$
then the theorem is proved. Otherwise, there exists a $f^{2}$-invariant leaf $I$ of $\mathcal{F}^{u}$. By theorem
\ref{TA} of section \ref{DF} we have that $S\cap I=\{p\}$. Thus,
$f^{2}$ has a fixed point $p$. Then $f$ has a periodic point and this contradicts infinitely expansiveness.
If $W^{u}$ does not separate $f^{-1}(W^{u})$ and $f(W^{u})$, applying lemma \ref{LP2} the existence of a $f^{2}$-invariant leaf $I$ of $\mathcal{F}^{u}$ is guaranteed. So we arrive to the same contradiction.

\item $W^{s}$ does not separate $f^{-1}(W^{s})$ and $f(W^{s})$. Applying lemma \ref{LP2} the existence of a $f^{2}$-invariant leaf $S$ of $\mathcal{F}^{s}$ is guaranteed. Then we can work with the same arguments used in the previous case.
\end{itemize}
Reciprocally, let´s consider the translation of the plane $T$ defined by $T(x,y)=(x+1,y)$. It was proved at the introduction that $T$ admits a Lyapunov metric function $D$ with the conditions required in this theorem. Now, let´s see the case when $f$ is conjugated to a translation $T$. Let us define a Lyapunov function for $f$ such as
$$L(p_{1},p_{2})=D(H(p_{1}),H(p_{2})),$$ where $D$ is the previous
defined Lyapunov metric function for $T$ and $H$ is a
homeomorphism from $\R^{2}$ over $\R^{2}$. It follows easily that
$L$ is a Lyapunov metric function for $f$ such that property {\bf HP} holds. This conclude the prove of the theorem.\ep

Using the results of \cite{G} and the last theorem we can state the following general characterization theorem:

\begin{teo}
Let $f$ be a homeomorphism of $\R^{2}$
which admits a Lyapunov function
$U:\R^{2}\times\R^{2}\rightarrow\R$ such that the following
properties hold:
\begin{itemize}
\item $U$ is a metric in $\R^{2} $ and induces the same topology
in the plane as the usual metric, \item For each point $x\in
\R^{2}$ and for each $k>0$ there exist points $y$ and $z$ on the
boundary of $B_{k}(x)$ such that $V(x,y)=U(f(x),f(y))-U(x,y)>0$ and
$V(x,z)=U(f(x),f(z))-U(x,z)<0$, respectively, \item $f$ has no
singularities.
\end{itemize}
Then,
\begin{itemize}
\item If $f$ admits a fixed point, then $f$ is conjugated to a linear hyperbolic automorphism
if and only if $U$ admits condition {\bf HP}; \item If $f$ is fixed point free then $f$ is topologically conjugate to a translation of the plane if and only if $U$ admits condition {\bf HP}
\end{itemize}
\end{teo}

\end{document}